\documentclass[12pt]{amsart}

\usepackage{amssymb}

\newtheorem{Thm}{Theorem}[section]

\newtheorem{Cor}[Thm]{Corollary}

\newtheorem{Lem}[Thm]{Lemma}

\newtheorem{Prop}[Thm]{Proposition}

\theoremstyle{definition}

\theoremstyle{remark}

\newcommand{\cobar}{\overline{\rm{co}}\,}

\def\ldots{\mathinner{\ldotp\ldotp\ldotp}}

\def\ldots{\mathinner{\cdotp\cdotp\cdotp}}

\def \cal{\mathcal}

\def \Bbb{\mathbb}

\def \diam{\text{diam }}

\begin{document}
\title{Szlenk indices and uniform homeomorphisms}
\author{G. Godefroy}
\address{
Equipe d'Analyse\\
Universit\'e Paris VI, Boite 186 \\
4, Place Jussieu\\
75252 Paris Cedex 05  \\
France}
\email[G. Godefroy]{gig@ccr.jussieu.fr}
\author{N. J. Kalton}
\address{Department of Mathematics \\
University of Missouri-Columbia \\
Columbia, MO 65211
}
\thanks{The second author was  supported by NSF grant DMS-9870027}

\email[N. J. Kalton]{nigel@math.missouri.edu}

\author {G. Lancien}
\address{Equipe de Math\'ematiques - UMR 6623, Universit\'e de
Franche-Comt\'e, F-25030 Besan\c con cedex}

\email{GLancien@vega.univ-fcomte.fr}


\begin{abstract} We prove some rather precise renorming theorems for
Banach spaces with Szlenk index $\omega_0.$  We use these theorems to
show the invariance of certain quantitative Szlenk-type indices under
uniform homeomorphisms.
\end{abstract}

\maketitle

\section{Introduction} \label{intro}

Classical results about super-reflexive Banach spaces include
Enflo-Pisier's renorming theorem
(\cite{E},\cite{P}) and Heinrich-Mankiewicz theorems on uniform and
Lipschitz homeomorphisms
\cite{HM}, from which it follows in particular that the moduli of uniform
convexity or uniform
smoothness of super-reflexive spaces is an invariant for uniform
homeomorphisms. This work is an
attempt to obtain similar results in the frame of non super-reflexive
spaces. We will in particular
relate, in a quantitative way, the Szlenk index with the existence of
equivalent $UKK^*$-renormings of the
space. These results extend and improve results in \cite{KOS}. We will
also show that the quantitative dependence on $\epsilon$
of the  Szlenk index (when finite) is an invariant under uniform
homeomorphisms.

We now turn to a detailed description of our results. In Section 2, the
{\it convex Szlenk index}
Cz is introduced and compared with the usual Szlenk index Sz. Precise
duality formulas, somewhat
related to duality between Orlicz spaces, are established which relate the
``$c_0$-like" behavior
of a Banach space
$X$ with the ``$l^1$-like" behavior of its dual $X^*$ (Proposition
\ref{orlicz}). Trees and
tree-maps are introduced in Section 3 as a tool to translate estimates on
the Szlenk index into
geometrical language. Section 4 is devoted to renormings.  It was
recently shown
\cite{KOS} that if Sz$(X)\le \omega_0$ then $X$ has an equivalent
$UKK^*$-renorming of some power type.  We prove more precise results.
 Trees are an
operative tool in the proof of the main renorming
theorem (Theorem \ref{renorm}) which provides the optimal relation between
the behavior of the convex Szlenk index
for any given value of $\tau>0$ and the existence of a 2-equivalent norm
with the $UKK^*$
property for this value of $\tau$.  In Theorem \ref{renorm2} we improve
the result of \cite{KOS} mentioned above by giving a precise bound on the
power-type.
 Note however that there is a small loss on the exponent in Theorem
\ref{renorm2}.
Such a loss is
unavoidable, as shown by the reflexive example from \cite{KOS}. This also
shows
Theorem \ref{renorm} cannot be improved to give a simultaneous renorming
which works for all $\tau>0.$
Summability of the Szlenk index is shown (Theorem \ref{sumszlenk}) to be
equivalent to the existence of a constant $K$ with
$\tau.Cz(X,\tau)<K$ for
any $\tau\in (0,1)$. Note that Corollary \ref{power} asserts that the
indices Cz and Sz have the same power-type behavior, at least
for $p>1$. Section 5 presents applications of these results to uniform
homeomorphisms. The main result of this section (Theorem
\ref{uniform}), whose proof uses the Gorelik principle ( Proposition
\ref{Gorelik}), asserts in a quantitative way that the
existence of equivalent
$UKK^*$ norms is invariant under uniform homeomorphisms. The main
application of this result (Theorem \ref{czinv}) is that the
convex Szlenk index is quantitatively invariant under uniform
homeomorphisms. This invariance is naturally obtained by a
combination of Theorem
\ref{uniform} with the renorming Theorem \ref{renorm}. Note that, although
the class of  spaces with separable dual is not stable
under uniform homeomorphisms
\cite{R}, the class of spaces with ``very separable" dual (that is, of
spaces with Szlenk index $\omega_0$) is by the above stable
under uniform homeomorphisms. It follows also (Theorem \ref{c0}) that a
Banach space which is
uniformly homeomorphic to $c_0$ is an isomorphic predual of $l^1$ with
summable Szlenk index.
We do not know whether a predual of $l^1$ with summable Szlenk index is
isomorphic to $c_0$.
Our last application (Theorem \ref{ellp}) concerns quotients of $l^p$ for
$p\in (2,\infty)$.

Some of the results of this work have been announced in \cite{GKL2}.

Acknowledgement: This work was initiated when the first and last named
authors were visiting the University of
Missouri-Columbia in 1997, and was concluded when the second named author
was visiting the Universit\'e de Besan\c con
in 1999. They are very grateful to these Institutions for their hospitality
and support.

\section{The Szlenk index and properties of norms}\label{index}

We first recall the definition of the Szlenk index and the Szlenk
derivation.
Suppose $X$ is a separable infinite-dimensional Banach space  and
$K$ is a weak$^*$-compact subset of
$X^*$.  If $\epsilon>0$ we let $\cal V$ be the set of all weak$^*$-open
subsets $V$ of $X^*$ such that $\diam V\cap K\le \epsilon$ and the
$\epsilon-$interior
$\iota_{\epsilon}K=K\setminus \cup\{V:V\in\cal V\}.$  (The set
$\iota_{\epsilon}K$ is often denoted $K_{\epsilon}'$ as in \cite{L},
\cite{L2}). We then define $\iota_{\epsilon}^{\alpha}K$ for any ordinal
$\alpha$ by
$\iota^{\alpha+1}_{\epsilon}K=\iota_{\epsilon}\iota_{\epsilon}^{\alpha}K$
and
$\iota^{\alpha}_{\epsilon}K=\cap_{\beta<\alpha}\iota_{\epsilon}^{\beta}K$
if $\alpha$ is a limit ordinal.

We denote by $B_{X^*}$ the closed unit ball of $X^*$. We then define
$\text{Sz}(X,\epsilon)$
(or
$\text{Sz}(\epsilon)$ if
no confusion can arise) to be the least countable ordinal $\alpha$
so that
$\iota_{\epsilon}^{\alpha}B_{X^*}=\emptyset,$ if such an ordinal exists.
Otherwise we will put $\text{Sz}(X,\epsilon)=\omega_1.$ The
{\it Szlenk index}  is defined by
$\text{Sz}(X)=\sup_{\epsilon>0}\text{Sz}(X,\epsilon).$ We recall that
$\text{Sz}(X)<\omega_1$ if and only if $X^*$ is separable.

Note that $\text{Sz}(X,\epsilon)\ge \epsilon^{-1}$ if $\epsilon>0$, and
compactness requires that $\text{Sz}(X,\epsilon)$ is not a limit ordinal. Thus
$\text{Sz}(X)=\omega_0$ is equivalent to $\text{Sz}(X,\epsilon)<\omega_0$
for every
$\epsilon>0,$ where $\omega_0$ denotes the first limit ordinal.

We  also introduce an alternative {\it convex Szlenk
index}. If
$K$ is compact and convex we may define
$\hat\iota_{\epsilon}K=\cobar
\iota_{\epsilon}K.$  Then $\text{Cz}(X,\epsilon)$ and $\text{Cz}(X)$
are defined as before, using instead this derivation.  Obviously
$\text{Cz}(X,\epsilon)\ge \text{Sz}(X,\epsilon).$  On the hand,
$\text{Cz}(X)<\omega_1$ is equivalent to  the separability of $X^*$; this
follows easily from the weak$^*$-dentability of weak$^*$-compact sets in
separable duals.

Henceforward we will only be interested in cases when
$\text{Sz}(X,\epsilon)$ and $\text{Cz}(X,\epsilon)$ are finite.  It will
therefore be useful to adopt the convention that
$\text{Sz}(X,\epsilon)=\infty$ if $\text{Sz}(X,\epsilon)\ge \omega_0$ and
similar convention for $\text{Cz}(X,\epsilon).$

Following \cite{KOS}, we will say that $X$ admits a {\it summable Szlenk
index} if there exists
a constant $K$ so that $\sum_{i=1}^n\epsilon_i \le K$ whenever
$\iota_{\epsilon_1}\ldots\iota_{\epsilon_n}B_{X^*}\neq \emptyset.$

The following lemma is proved in \cite{L2}, p.57, or \cite{L}:

\begin{Lem}\label{submult} If $0<\epsilon,\eta\le 1,$
 then
$\text{Sz}(X,\epsilon\eta)\le \text{Sz}(X,\epsilon)\text{Sz}(X,\eta).$\end{Lem}

Note that this implies that $\text{Sz}(X)=\omega_0$ if and only if
$\text{Sz}( X,\epsilon)<\infty$ for any fixed $0<\epsilon<1.$

Another immediate consequence of this is that if $\text{Sz}(X)=\omega_0$
then

$$\lim_{\epsilon\to 0}\frac{\log \text{Sz}(X,\epsilon)}{|\log \epsilon|
}=p$$
exists where $1\le p<\infty$.  We will call $p=p_X$ the {\it Szlenk power
type} of $X.$  We also have that if $\delta>0$ then $\text{Sz}(X,\epsilon)\le
C\epsilon^{-p-\delta}$ for some suitable constant $C$ so that
$\text{Sz}(X,\epsilon)$ grows at a power rate.  In fact we can also
define
$p_X$
as the infimum of all $q$ so that $\epsilon^q\text{Sz}(X,\epsilon)$ is
bounded.

Next we note

\begin{Lem}\label{convex} If $0<\epsilon<1$ and $n\in\Bbb N$ are such
that $n\epsilon\le 1$ then:\newline (1)
$\text{Sz}(X,\epsilon)-1\ge
n(\text{Sz}(X,n\epsilon)-1).$ \newline (2)
$\text{Cz}(X,\epsilon)-1\ge n(\text{Cz}(X,n\epsilon)-1).$\end{Lem}

\begin{proof} Note that if $j<\text{Sz}(X,n\epsilon)$ then for any $m\ge 2$,
$$(m-1)B_{X^*}+\iota_{n\epsilon}^jB_{X^*}\subset
\iota_{n\epsilon}^j(mB_{X^*})$$  Hence $(m-1)B_{X^*}\subset
\iota_{n\epsilon}^{l}mB_{X^*}$ where $l=\text{Sz}(X,n\epsilon)-1.$
Iterating we
obtain $\iota_{n\epsilon}^{nl}(nB_{X^*})\neq \emptyset$ and so
$\text{Sz}(X,\epsilon) \ge nl+1.$ This implies the result for the Szlenk index
and the convex version is similar.

\end{proof}

We will also need the following elementary fact:

\begin{Lem}\label{iso} Suppose $X$ and $Y$ are isomorphic Banach spaces.
  Then if $d=d(X,Y)$ we have
$$ \text{Sz}(X,d\epsilon)\le \text{Sz}(Y,\epsilon)\le
\text{Sz}(X,d^{-1}\epsilon)$$ and
$$ \text{Cz}(X,d\epsilon)\le \text{Cz}(Y,\epsilon)\le
\text{Cz}(X,d^{-1}\epsilon).$$

\end{Lem}

There is one important advantage of the convex Szlenk index which is
established in the following lemma. Note that part (2) shows that if the
convex Szlenk index is $O(\tau^{-1})$, then it is actually
summable.

\begin{Lem}\label{summable} Suppose $X$ is a separable Banach space,
$0\le \epsilon_k\le 1$ for $1\le k\le N$ and $s=\sum_{i=1}^N\epsilon_i.$
\newline
(1) If
$\iota_{\epsilon_1}\ldots\iota_{\epsilon_N}B_{X^*}\neq
\emptyset$ and $0<\tau\le
\frac{s}{2N}$
$$ \sum_{\substack{k\ge 0\\ 2^k\tau\le 1}}2^k(\text{Sz}(X,2^k\tau)-1)\ge
\frac{s}{2\tau}.$$\newline
(2) If $\hat\iota_{\epsilon_1}\ldots\hat\iota_{\epsilon_N}B_{X^*}\neq
\emptyset$ and
 $0<\tau\le
s/(4N)$,

$$\text{Cz}(X,\tau))-1\ge \frac{s}{8\tau}.$$

\end{Lem}

\begin{proof}

For the first part, let us write $\alpha_i\le \epsilon_i\le \beta_i$
where $\alpha_i=2^{k_i}\tau$ and $\beta=2^{k_i+1}\tau$ with $k_i$ a
nonnegative integer, or $\alpha_i=0$ and $\beta_i=\tau.$ Then
$\beta_i=\tau+ \sum_{0\le j\le k_i}2^j\tau.$  Now observe that
$|\{j:\epsilon_j\ge
2^k\tau\}|<
\text{Sz}(X,2^k\tau)$ and so
$$ s\le N\tau+
 \sum_{\substack{k\ge 0\\ 2^k\tau\le
1}}2^k\tau(\text{Sz}(X,2^k\tau)-1).$$
This immediately gives the conclusion of (1).

 For the second part, first suppose $C$ is a weak$^*$-compact and convex
subset
of $X^*$.  We will argue that if $\epsilon>0$ and $k\in\Bbb N$ then
$\hat\iota_{2k\epsilon}C\subset\hat \iota_{\epsilon}^kC.$  Indeed suppose
$x^*\in\iota_{2k\epsilon}C.$  Then there is a sequence $x_n^*\in C$ with
$\|x_n^*-x^*\|\ge k\epsilon$ and $x_n^*\to x^*$ weak$^*$.   Now for any
$n_1<n_2<n_3\ldots<n_k$ we have $\frac1k(x^*_{n_1}+\cdots+x^*_{n_k})\in C.$
  Letting $n_k\to\infty$ we obtain that
$\frac1k(x^*+x^*_{n_1}+\cdots+x^*_{n_{k-1}})\in \iota_{\epsilon}C.$
Letting $n_{k-1}\to\infty$ and repeating we obtain that $x^*\in
\iota_{\epsilon}^kC.$
In particular we observe that
$\hat\iota_{2k\epsilon}C\subset\hat\iota_{\epsilon}^{k}C.$

Now suppose $\epsilon_1,\epsilon_2,\ldots,\epsilon_N>0$ and
$\sum_{k=1}^N\epsilon_k=s.$  Suppose $\tau\le s/4N$ and let
$m_k=[\epsilon_k/2\tau].$  Then
$$\hat\iota_{\epsilon_1}\ldots\hat\iota_{\epsilon_N}B_{X^*} \subset
\hat\iota_{\tau}^{m_1+\cdots+m_N}B_{X^*}.$$

Now $$m_1+\cdots+m_N= \sum_{\epsilon_k\ge 2\tau}\left[\frac
{\epsilon_k}{2\tau}\right]\ge
\frac1{4\tau}
\sum_{\epsilon_k\ge
2\tau}\epsilon_k.$$

But $\tau\le \frac{s}{4N}$ yields $\sum_{\epsilon_k\le 2\tau}\epsilon_k\le
\frac{s}{2}.$  Hence $m_1+\cdots+m_N\ge \frac{s}{8\tau}.$

\end{proof}

Let us now turn to renormings.  We need the following elementary lemma:

\begin{Lem}\label{weaknull} Suppose $X$ is an infinite-dimensional Banach
space with separable
dual and that
$(x_n^*)$ is a sequence in $X^*$ with $\lim_{n\to\infty}x_n^*=0$
weak$^*$.  Then there is a sequence $(x_n)$ in $X$ with $\|x_n\|\le 1,$
$\lim_{n\to \infty}x_n=0$ weakly and
$$ \liminf\langle x_n,x_n^*\rangle \ge \frac12\liminf \|x_n^*\|.$$

\end{Lem}

\begin{proof}

We observe that by Lemma 2.3 of \cite{KW} we have
$$\liminf_{n\to\infty}d(x_n^*,F)\ge
\frac12\liminf_{n\to\infty}\|x_n^*\|$$ if
$F$ is any
finite-dimensional subspace of $X^*$.  Hence, since $X^*$ is separable we
find an increasing sequence (not necessarily strictly increasing) of
finite-dimensional subspaces $(F_n)$ so that $\cup F_n$ is dense in $X^*$
and $\liminf d(x_n^*,F_n)\ge \frac12\liminf_{n\to\infty}\|x_n^*\|.$

There exist
$x_n\in X$ with $\|x_n\|=1,\ x_n\in F_n^{\perp}$ and
$\langle x_n,x^*_n\rangle>d(x_n^*,F_n)-\frac1n$ and this gives the
conclusion.

\end{proof}

\begin{Prop}\label{renorming}.  Suppose $X$ is a separable Banach space
 and
 $0<\sigma,\tau < 1.$  Consider the
following statements:

\begin{enumerate}

\item
If $x^*\in X^*$ with $\|x^*\|=1$, $\lim x_n^*=0$ weak$^*$ with
$\lim_{n\to\infty}\|x_n^*\|=\tau$ then
$$ \liminf_{n\to\infty}\|x^*+x_n^*\|\ge 1+\sigma\tau.$$
\item  $X^*$ is separable and
if
$x\in X$ with
$\|x\|=1$,
$\lim x_n=0$ weakly  with
$\lim\|x_n\|=\sigma$ then

$$ \limsup_{n\to\infty}\|x+x_n\|\le 1+\sigma\tau.$$

\item $X^*$ is separable and
if $x^*\in X^*$ with $\|x^*\|=1$, $\lim x_n^*=0$ weak$^*$ with
$\lim_{n\to\infty}\|x_n^*\|=6\tau$ then
$$ \liminf_{n\to\infty}\|x^*+x_n^*\|\ge 1+\sigma\tau.$$
\end{enumerate}

Then $(1)\Rightarrow (2)\Rightarrow (3)$.
 Furthermore
(1) implies that
$$\text{Sz}(X,2\tau)\le \text{Cz}(X,2\tau)\le \sigma^{-1}\tau^{-1}+1.$$
(In particular if
$2\tau<1$ then
$\text{Sz}(X)=\omega_0.$) \end{Prop}

\begin{proof}

Note  that (1) implies that $\iota_{2\tau}B_{X^*} \subset
(1-\sigma\tau)B_{X^*}$ and this immediately yields $\text{Cz}(X,2\tau)\le
\sigma^{-1}\tau^{-1}+1.$ Then  the last statement follows
from Lemma \ref{submult}.

First assume (1) holds
and
$x,x_n$ are chosen as in
(2).
It is enough to show it is impossible that $\lim_{n\to\infty}\|x+x_n\|>
1+\sigma\tau.$  Suppose this holds.  Then we can pick $y_n^*\in
B_{X^*}$ with $\lim_{n\to\infty}\langle x+x_n,y_n^*\rangle
>1+\sigma\tau.$  Passing to a subsequence we can suppose $y_n^*$
converges weak$^*$ to some $x^*\in B_{X^*}$ and then put
$x_n^*=y_n^*-x^*.$  We can assume that
$\lim_{n\to\infty}\|x_n^*\|=\theta$ exists.

If $\theta\le \tau$ then
$$\lim_{n\to\infty}\langle x+x_n,x^*+x_n^*\rangle \le \langle
x,x^*\rangle +\sigma\tau \le 1+\sigma\tau.$$
On the other hand if $\theta>\tau$ then we have, using the convexity of
the norm,

$$ 1=\liminf_{n\to\infty}\|x^*+x_n^*\| \ge \|x^*\| +\sigma \theta.$$

Hence $\|x^*\|\le 1-\sigma\theta$ and so
$$\lim_{n\to\infty}\langle x+x_n,x^*+x_n^*\rangle \le
1-\sigma\theta+\sigma\theta=1.$$  This gives us the required
contradiction.

 Assume  (2) holds.
 Suppose that
$x^*,x_n^*$ are chosen as in
(3).  It will be enough to show the conclusion for some subsequence.

Then, given $\epsilon>0$ we can choose $x\in X$ with $\|x\|=1$ and
$x^*(x)>1-\epsilon.$  Now by Lemma \ref{weaknull}, we can,
by passing to a subsequence, assume there exist
a weakly null sequence $x_n\in X$ with $\|x_n\|=1$ and $\liminf x_n^*(x_n)\ge
3\tau.$  Now
$$ \liminf \langle x+\sigma x_n,x^*+x_n^*\rangle \ge
1-\epsilon+3\sigma\tau$$ and so
$$ 1-\epsilon+3\sigma\tau\le (1+\sigma\tau)
\liminf\|x^*+x_n^*\|.$$ Hence
 letting
$\epsilon\to 0$ we have

$$ \liminf \|x^*+x_n^*\| \ge \frac{1+3\sigma\tau}{1+\sigma\tau}\ge
1+\sigma\tau.$$

\end{proof}

\begin{Prop}\label{ell} Let $X$ be a separable Banach space not
containing a copy of $\ell_1,$ and suppose $0<\sigma<1$ and
$0<\tau<\frac12.$ Suppose  that whenever
$x\in X$ with $\|x\|=1$ and $\lim x_n=0$ weakly with
$\lim_{n\to\infty}\|x_n\|=\sigma$ then

$$ \limsup_{n\to\infty}\|x+x_n\|\le 1+\sigma\tau.$$

Then $X^*$ contains no proper norming subspaces and hence is
separable.\end{Prop}

\begin{proof} Suppose $x^{**}\in X^{**}$ is of norm one and
$\{x^*\in X^*:x^{**}(x^*)=0\}$ is a norming subspace of $X^*.$  Then by
the Odell-Rosenthal theorem \cite{OR} there is a weakly Cauchy sequence
$(x_n)$ in $X$ with $\|x_n\|=1$ and $x_n\to x^{**}$ weak$^*$.

 If $m\in\Bbb N$ then, for any $\epsilon>0,$ and for each
$n>m$ we can choose
$e^*\in X^*$ with $\|e^*\|=1$, $x^{**}(e^*)=0$ and
$e^*(x_m+\frac12\sigma x_n)>\|x_m+\frac12\sigma x_n\|-\epsilon.$  Then we
can find
$k(n)>n$ so that
$e^*(\sigma x_{k(n)})<\epsilon.$  We conclude that

$$ \|x_m+\frac12\sigma(x_n-x_{k(n)})\| >\|x_m+\frac12\sigma
x_n\|-2\epsilon.$$

Letting $n\to\infty$, since $(x_n-x_{k(n)})\to 0$ weakly we have
$$ \limsup_{n\to\infty}\|x_m+\frac12\sigma x_n\| \le
1+\sigma\tau+2\epsilon.$$

Thus

$$ \|x_m+\frac12\sigma x^{**}\|\le 1+ \sigma\tau.$$

Then $m\to\infty$ we have $1+\frac12\sigma\le 1+\sigma\tau$ which is a
contradiction.\end{proof}

At this point we introduce some terminology.  Let $f,g$ be continuous
monotone increasing functions on $[0,1]$ which satisfy $f(0)=g(0)=0.$
We will say that $f$ {\it C-dominates} $g$ if
$f(\tau)\geq g(\tau/C)$ for every $0\le \tau\le 1.$  We will say that $f,g$
are {\it
C-equivalent} if $f$ $C$-dominates $g$ and $g$ $C$-dominates $f.$

For any such monotone increasing function $f$ we denote by $f^*$ its
dual Young's function i.e.
$$ f^*(s)=\sup\{st-f(t):\ 0\le t\le 1\}.$$

Notice that if $f$ $C$-dominates $g$ then $g^*$ $C$-dominates $f^*.$
Note also that $f^*$ is a convex function.

Now if $X$ is any separable Banach space we define for
$0\le\sigma\le 1,$ $\rho(\sigma)=\rho_X(\sigma)$ to be the least constant
so that

$$ \limsup_{n\to\infty}\|x+x_n\|\le 1+\rho_X(\sigma)$$ whenever
$\|x\|=1,$ $\lim x_n=0$ weakly  and $\limsup_{n\to\infty}\|x_n\|\le
\sigma.$  We define $\theta(\tau)=\theta_X(\tau)$ for $0\le\tau\le 1$ to
be the greatest constant so that
$$ \liminf_{n\to\infty}\|x^*+x_n^*\|\ge 1+\theta_X(\tau)$$ whenever
$x^*,x_n^*\in X^*,\ \|x^*\|=1,\ \lim x_n^*=0$ weak$^*$ and
$\liminf_{n\to\infty}\|x_n^*\|\ge\tau.$

We then define $\varphi(\sigma)$ by
$\varphi(\sigma)=\inf\{\rho_Y(\sigma):\ d(X,Y)\le 2\}$ and
$\psi(\tau)=\sup\{\theta_Y(\tau):\ d(X,Y)\le 2\}.$

We can now summarize Proposition \ref{renorming} and
Proposition \ref{ell}:

\begin{Prop}\label{orlicz} Let $X$ be a separable Banach space not
containing $\ell_1.$ Then:
\begin{enumerate}\item Each of the functions $\rho(t)/t$, $\theta(t)/t,$
$\varphi(t)/t$ and $\psi(t)/t$ is monotone increasing on $(0,1).$
\item $\theta$ is 4-equivalent to $\rho^*$ and $\psi$ is
4-equivalent to $\varphi^*.$
\item  $\rho$ is 8-equivalent to $\theta^*$ and
$\varphi$ is 8-equivalent to $\psi^*$
\item If $0<t\le 1$ then $$\text{Cz}(X,t)\le 1+
\frac{1}{\psi(t/4)}.$$
\end{enumerate}\end{Prop}

\begin{proof} (1) follows trivially from convexity considerations.

(2) It follows from Proposition \ref{ell} that if $X^*$ is not separable,
then $\rho^*(\tau)=0$
for all $0\leq \tau \leq \frac 12$ and then it is obvious that $\theta$
2-dominates
$\rho^*$. If $X^*$ is separable, we pick
$(x_{n}^*)$ in
$X^*$ such that $\liminf_{n\to\infty}\|x_n^*\|\ge\tau$ and $$
\liminf_{n\to\infty}\|x^*+x_n^*\|=1+\theta_X(\tau)$$
By Lemma \ref{weaknull} there is a weakly null sequence $(x_n)$ with
$\limsup_{n\to\infty}\|x_n\|\le
\sigma$ and $\lim x_{n}^*(x_n)=\sigma\tau/2$. It follows easily
 that
$$1+\sigma\tau/2
\le
(1+\rho(\sigma))(1+\theta(\tau))$$ so that
$$ \sigma\tau/2\le \rho(\sigma)+\theta(\tau)+\rho(\sigma)\theta(\tau).$$

Now since $\theta(\tau)\le 1$ we obtain
$$ \theta(\tau)\ge \sigma\tau/2-2\rho(\sigma)\ge
2\rho^*(\tau/4).$$  Hence $\theta$ 4-dominates $\rho^*.$

 The same considerations show that $\rho$ 4-dominates $\theta^*$,
$\varphi$ 4-dominates $\psi^*$ and $\psi$ 4-dominates  $
\varphi^*.$

Now if $0<\tau<1$ pick $\sigma=2\theta(\tau/2)/\tau.$  Then by
Proposition \ref{renorming} we have $\rho(\sigma) \le \frac12\sigma\tau.$
Hence
$$\theta(\tau/2)=\frac12\sigma\tau\le
\sigma\tau-\rho(\sigma)\le \rho^*(\tau).$$  Thus $\rho^*$ 2-dominates
$\theta.$ The proof for $\psi$ in place of $\theta$ and $\varphi$ in
place of
$\rho$ is similar.

(3) We deduce from (2) that $\theta^*$ is 4-equivalent to $\rho^{**}$.
Next let $$\tilde\rho(t)=\int_0^t\frac{\rho(\tau)}{\tau}d\tau.$$
Then $\rho(t/2)\le\tilde\rho(t)\le\rho(t).$  Hence since $\tilde\rho$ is
convex we have
$\rho^{**}(t)\ge \tilde\rho(t)\ge \rho(t/2).$  Hence $\rho$ is
2-equivalent to $\rho^{**}.$  The argument for
$\varphi$ and $\psi$ is similar.

(4) This is an immediate deduction from Proposition \ref{renorming} and
Lemma \ref{iso}.
\end{proof}

\section{Trees and tree-maps}\label{trees}

Consider the set ${\cal F}{\Bbb N}$ of all finite subsets of $\Bbb N$
with the following partial order.  If $a=\{n_1,n_2,\ldots,n_k\}$ where
$n_1<n_2<\ldots<n_k$ and $b=\{m_1,m_2\ldots,m_l\}$ where
$m_1<m_2<\cdots<m_l,$ then $a\le b$ if and only if $k\le l$ and
$m_i=n_i$ where $1\le i\le k$ (i.e. $a$ is an initial segment of $b.$)
We say that $b$ is a {\it successor} of $a$ if $|b|=|a|+1$ and $a\le b;$
the collection of successors of $a$ is denoted by $a+$.
If $a\neq\emptyset$ then $a-$ denotes the unique predecessor of $a$ i.e.
$a$ is a successor of $a-.$ Let
$S$ be a subset of
$\cal F\Bbb N.$ We will say that
$S$ is a
{\it full tree}
if we have
 \begin{enumerate} \item $\emptyset \in S.$
  \item  Each
$a\in S$ has infinitely many
successors in $S.$
\item If $a\in S$ and $\emptyset\neq a\in S$ then
$a-\in S$.
\end{enumerate}

  It is easy to see that any full tree
 is isomorphic as an ordered set to $\cal F\Bbb N.$
 If $S$ is any full tree we will say that a sequence
$\beta=\{a_n\}_{n=0}^{\infty}$ is a {\it branch} of $S$ if $a_n\in S$
for all
$n,$ $a_0=\emptyset$ and $a_{n+1}$ is a successor of $a_n$ for all
$n\ge
0.$

Now let $V$ be a vector space.  We define a {\it tree-map} to be a map
$a\mapsto x_a$ defined on a full tree $S$ with the properties that
$x_{\emptyset}=0$ and for every branch $\beta$ the set
$\{a:\ x_{a}\neq 0:\ a\in\beta\}$ is finite.  Given any tree-map we
define a {\it height function} $h$ which assigns to each $a$ a countable
ordinal; to do this we define $h(a)=0$ if $x_{b}=0$ for $b\ge a$ and
then inductively $h(a)$ is defined by $h(a)\le \eta$ if and only if
$h(b)<\eta$ for every $b>a.$  The {\it height} of the tree-map is
defined to be
$h(\emptyset).$  Note that the tree-map $a\mapsto x_{a}$ has finite
height
$m\le n$ if and only if $x_{a}=0$ whenever $|a|>n.$ For a recent work on
trees in Banach spaces we refer to \cite{AJO}.

The following easy lemma is a restatement of the fact that certain types
of games (which are not used in this paper) are determined.

\begin{Lem}\label{determine}Suppose $(x_a)_{a\in S}$ is a tree-map and
that $A$ is any subset of $V.$  Then either there is a full tree
$T\subset S$ so that $\sum_{a\in\beta}x_a\in A$ for every branch
$\beta\subset T$ or there is a full tree $T\subset S$ so that
$\sum_{a\in\beta}x_a\notin A$ for every branch $\beta\subset
T.$\end{Lem}

\begin{proof} For each countable ordinal $\eta$ we define a subset
$B_{\eta}$ of
$\{a\in S: h(a)=\eta\}$ as follows. If $\eta=0$ let $a\in B_0$ if
 $\sum_{b\le a}x_a\in A.$  Then inductively if
$h(a)=\eta$ we say $a\in B_{\eta}$ if $a$ has infinitely many
successors
$b$ with $b\in B_{h(b)}.$  Let $B=\cup_{\eta}B_{\eta}.$  If $\emptyset\in
B$ then an easy induction argument produces a full tree $T\subset S$ with
$\sum_{a\in\beta} x_{a}\in A$ for every $\beta\subset T.$  Otherwise the
set $T=S\setminus B$ is a full tree  with the property that $\sum_{a\in
\beta}x_{a}\notin A$ for every $\beta\subset T.$\end{proof}

We now consider tree-maps with values in a Banach space $X.$

\begin{Lem}\label{reduce} Suppose $(x_a)_{a\in S}$ is a bounded tree-map
in $X$ of finite height $n.$  Then, given $\delta>0$ we can find a full
tree $T\subset S$ and $\epsilon_1,\ldots,\epsilon_n\ge 0$ so that if
$a\in T$ and $|a|=k\le n$ then $\epsilon_k\le \|x_a\|\le
\epsilon_k+\delta.$\end{Lem}

\begin{proof} One easy way to prove this is to consider $V=\Bbb R^n$ with
canonical basis $e_1,\ldots,e_n$ and then the tree-map
$ u_a=
\|x_a\|e_k$ if $|a|=k\le n$ and $u_a=0$ if $|a|>n.$  The lemma follows
from the Heine-Borel theorem and repeated applications of Lemma
\ref{determine}. \end{proof}

 If $\tau$
is a topology on $X$ (e.g. the weak topology or for dual spaces the
weak$^*$-topology) we say that a tree-map $(x_a)_{a\in S}$ is {\it
$\tau-$null} if for every $a\in S$ the set $\{x_b\}_{b\in a+}$ is a
$\tau$-null sequence.

\begin{Lem}\label{dual} Suppose $X$ is a Banach space and that
$(x_a)_{a\in S}$
is  a weakly null tree-map in $X$ and $(x_a^*)_{a\in S}$ is a
weak$^*$-null
tree map in $X^*$.  Then for any $\delta>0$ there is a full tree
$T\subset S$ so that for any branch $\beta\subset T$ we have
$$ |\langle \sum_{a\in\beta}x_a,\sum_{a\in\beta}x^*_a\rangle
-\sum_{a\in\beta}\langle x_a,x_a^*\rangle|\le \delta.$$
\end{Lem}

\begin{proof} Let $T$ be the set of $a$ so that if $b< c\le a$ then
$|\langle x_b,x^*_c\rangle|,|\langle x_c,x^*_b\rangle|\le
|c|^{-1}2^{-|c|-1}\delta.$  It is not difficult to see that this is
a full tree and that if $\beta$ is a branch in this tree then
$$ |\langle \sum_{a\in\beta}x_a,\sum_{a\in\beta}x^*_a\rangle
-\sum_{a\in\beta}\langle x_a,x_a^*\rangle|= |\sum_{b<
a}\sum_{a\in\beta}\langle x_a,x_b^*\rangle+\langle
x_b,x_a^*\rangle| \le \delta.$$
                 \end{proof}

\begin{Prop}\label{Szlenk} Let $X$ be a separable Banach space. In order
that
$$\iota_{\epsilon_1}\iota_{\epsilon_2}\ldots\iota_{\epsilon_n}B_{X^*}\neq
\emptyset
$$ it
 is necessary  that there
exists
a weak$^*$-null tree-map $(x^*_a)_{a\in S}$ in $X^*$ with $\|x^*_a\|\ge
\frac14\epsilon_{|a|}$
for $1\le |a|\le n$ and so that $\|\sum_{a\in\beta}x^*_a\|\le 1$ for
every branch $\beta,$ and sufficient that
  there
exists
a weak$^*$-null tree-map $(x^*_a)_{a\in S}$ in $X^*$ with $\|x^*_a\|\ge
\epsilon_{|a|}$
for $1\le |a|\le n$ and so that $\|\sum_{a\in\beta}x^*_a\|\le 1$ for
every branch $\beta$.
\end{Prop}

\begin{proof}
First assume
$\iota_{\epsilon_1}\iota_{\epsilon_2}\ldots\iota_{\epsilon_n}B_{X^*}=K\neq
\emptyset.$  Then
$$ 0\in
\frac12K+\frac12B_{X^*}\subset\iota_{\epsilon_1/2}\ldots\iota_{\epsilon_n/2}B_{X
^*}$$
Now there exists a sequence $x^*_k$
converging to $0$ weak$^*$ so that $\|x_k^*\|\ge\frac14\epsilon_1$ and
$x_k^*\in
\iota_{\epsilon_2/2}\iota_{\epsilon_3/2}\ldots\iota_{\epsilon_{n}/2}
B_{X^*}$ since otherwise
there is a weak$^*$-open neighborhood of $0$
relative to
$\iota_{\epsilon_2/2}\iota_{\epsilon_3/2}\ldots\iota_{\epsilon_{n}/2}
B_{X^*}$
 of diameter
less than $\epsilon_1/2.$ Then for each
$k\in
\Bbb N$ we find a sequence $(x^*_{k_1,k_2})_{k_2>k_1}$ so that
$\lim_{k_2}x^*_{k_1,k_2}=0$ weak$^*$, $\|x^*_{k_1,k_2}\|\ge \frac14
\epsilon_{2}$ and $x^*_{k_1}+x^*_{k_1,k_2}\in
\iota_{\epsilon_3/2}\ldots\iota_{\epsilon_{n}/2}B_{X^*}.$  This procedure
can then
be iterated to define $x^*_a$ if $|a|\le n$.  Setting $x_a^*=0$ if
$|a|>n$ we obtain the desired tree-map.

The converse is equally easy.  Obviously if $(x^*_a)_{a\in S}$ is the
given tree-map then we have $\|\sum_{a\le b}x^*_a\|\le 1$ for all $b.$
It then follows inductively that $\sum_{b\le a}x^*_b\in
\iota_{\epsilon_{k+1}}\iota_{\epsilon_{k+2}}\ldots\iota_{\epsilon_n}
B_{X^*}$ whenever $|a|=k.$  Setting
$k=0$ gives the result.\end{proof}

\section{$UKK^*$-renormings}\label{UKK}

Suppose $X$ is a separable Banach space. If
$\sigma>0$ we define $N=N(\sigma)$ to be the least integer $N$ so that
there exists a weakly null tree-map $(x_a)_{a\in S}$ in $X$ of height
$N+1$ such that $\|x_a\|\le \sigma$ for every $a\in S$ and
$\|\sum_{a\in\beta}x_a\|> 1$ for every branch $\beta.$ (We put
$N(\sigma)=\infty$ if no such integer exists.)

 Notice that
$N(\sigma)>\sigma^{-1}-1.$

{\it Remark.} It follows from Lemma
\ref{determine} that if $k\le N(\sigma)$ then for every weakly null
tree-map
$(x_a)_{a\in S}$ of height $k$ and such that $\|x_a\|\le \sigma$ for all
$a\in S$ there is a full tree $T\subset S$ so that  $\|\sum_{a\in
\beta}x_a\|\le 1$ for every branch.

We refer to \cite{DGK} and references therein for the uniform Kadec-Klee
property and its dual version. Here we denote by $UKK^*$
what was named there  weak-star $UKK$. Also, we consider it as a property
of a Banach space which can be checked on its dual, rather
than a property of a dual space which refers to a given predual. Our main
reason for introducing the quantity
$N(\sigma)$ is to obtain a renorming of the space $X$ with an approximate
$UKK^*$-condition.

\begin{Thm}\label{UKK1}Suppose $X$ is a separable Banach space.  Then
for any $\sigma>0$ if $N(\sigma)<\infty$  there is a norm
$|\cdot|$ on
$X$ satisfying $\frac12\|x\|\le |x|\le \|x\|$ and
$$ \limsup_{n\to\infty}|x+x_n|\le 1+ \frac{1}{N(\sigma)}$$ whenever
$|x|=1$ and
$(x_n)$ is a sequence satisfying $\lim_{n\to\infty}|x_n|=\frac12\sigma$
and
$\lim_{n\to\infty}x_n=0$ weakly.\end{Thm}

\begin{proof}Define $f_0(x)=\|x\|$ and then for $k>0$ define $f_k(x)$ to
be the infimum of all $\lambda>0$ so that whenever $(x_a)_{a\in S}$ is a
weakly null tree-map of height $k$ with $\|x_a\|\le\sigma$ for all
$a\in S$ then there is a full subtree $T\subset S$
so that
$\|x+\sum_{a\in\beta}x_a\|\le \lambda$ for every branch $\beta.$
We observe first that $(f_k(x))_{k=0}^{\infty}$ is an increasing
sequence, and that $f_k(x)=f_k(-x).$ Next notice that
$$ |f_k(x)-f_k(y)|\le \|x-y\|$$ by an elementary calculation, which we
omit.
We also observe that $f_N(x)\le \|x\|+1$ where $N=N(\sigma).$

We next claim that $f_k$ is convex.
Indeed  let $u=t x+(1-t)y,$ where $0<t<1.$
Suppose
$\lambda>f_k(x)$ and $\mu>f_k(y).$  Let $(x_a)_{a\in S}$ be any weakly
null tree-map of height $k$ with $\|x_a\|\le \sigma$ for all $a\in S.$

Then we can find a full subtree $T_1\subset S$ so that for every branch
$\beta$ we have
$$ \|x+\sum_{a\in\beta}x_a\|\le \lambda$$
and then a full subtree $T_2\subset T_1$ so that for every branch
$\beta\subset T_2$
$$ \|y+\sum_{a\in\beta}x_a\| \le \mu.$$
Obviously for every branch $\beta\subset T_2$
$$ \|u+\sum_{a\in\beta}x_a\| \le t\lambda+(1-t)\mu$$
so that $f_k(u)\le t\lambda+(1-t)\mu.$

Next we note that if $\|x_n\|\le \sigma$ and $\lim_{n\to\infty}x_n=0$
weakly then $$\limsup f_k(x+x_n)\le f_{k+1}(x)$$ for all $k\ge 0.$  Indeed
for $k=0$ this is obvious.  If $k>0$ assume that
$\lambda <\limsup_{n\to\infty}f_k(x+x_n).$  By passing to a
subsequence we can suppose $\lambda<f_k(x+x_n)$ for every $n.$ Then for
each
$n$ there is a weakly null tree-map  $(y^{(n)}_a)_{a\in S_n}$ of height
$k$ so that $\|y^{(n)}_a\|\le \sigma$ for all $a\in S_n$ and $$
\|x+x_n+\sum_{a\in\beta}y^{(n)}_a\| >\lambda$$ for every branch
$\beta\subset S_n.$

Now let $T$ be the tree consisting of all sets $\{m_1,\ldots,m_l\}$
where
$m_1<m_2<\cdots<m_l$ such that if $l>1$ then $\{m_2,\ldots,m_l\}\in
S_{m_1}.$  We define a weakly null tree-map of height $k+1$ by
$$ z_{m_1,\ldots,m_l}=\begin{cases} x_{m_1}  \text{ if } l=1\\
y^{(m_1)}_{m_2,\ldots,m_l} \text{ if } l>1.\end{cases}$$
Then for every branch $\beta\subset T$  we have
$$ \|x+\sum_{a\in\beta}z_a\| >\lambda$$
so that $f_{k+1}(x)\ge \lambda$.
This implies our claim.

Now let us set $g(x)=\frac1N\sum_{k=0}^{N-1}f_k(x).$  Then $g$ is convex,
and $\|x\|\le g(x)\le \|x\|+1.$ Further
if $(x_n)$ is weakly null with $\|x_n\|\le\sigma$ for all $n$ then
\begin{equation}\begin{align*} \limsup_{n\to\infty}g(x+x_n) &\le \frac1N
\sum_{k=1}^Nf_k(x)\\ &\le g(x)+\frac1N(f_N(x)-f_0(x))\\ &\le g(x)
+\frac1N.\end{align*}\end{equation}

Let $|\cdot|$ be the Minkowski functional of the set $\{x:g(x)\le 2\}.$
Then it is clear that $\frac12\|x\|\le |x|\le \|x\|.$  Suppose $|x|=1$
and
$(x_n)$ is a weakly null sequence with $|x_n|\le \frac12\sigma.$  Then
$\|x_n\|\le \sigma.$  Thus
$$ \limsup_{n\to\infty} g(x+x_n)\le g(x)+\frac1N \le 2+\frac1N.$$
Now $g(0)< 1$ and so from the convexity of $g$ it follows that
$$ \limsup_{n\to\infty} g(\frac{N}{N+1}(x+x_n))< 2$$
whence
$$ \limsup_{n\to\infty}|x+x_n| \le 1+\frac1N.$$

\end{proof}

\begin{Thm}\label{UKK2}Suppose $X$ is a separable Banach space.
Suppose $\sigma>0$ is such that $N(\sigma)=\infty.$  Then, for any
$\epsilon>0$ there is a norm
$|\cdot|$ on
$X$ satisfying $\frac12\|x\|\le |x|\le \|x\|$ and
$$ \limsup_{n\to\infty}|x+x_n|\le 1+ \epsilon$$
whenever $|x|=1$ and
$(x_n)$ is a sequence satisfying $\lim_{n\to\infty}|x_n|=\frac12\sigma$
and
$\lim_{n\to\infty}x_n=0$ weakly.\end{Thm}

\begin{proof}The proof is almost identical to the proof of the  preceding
Theorem
\ref{UKK1}, except that one considers
$g_m(x)=\frac1m\sum_{k=0}^{m-1}f_k(x)$ for arbitrarily large choices of
$m.$  We omit the details. \end{proof}

Notice that the preceding two theorems allow us to say in the language of
Section \ref{index} that $N^{-1}$ 2-dominates
$\varphi$ where $N^{-1}(\sigma)=(N(\sigma))^{-1}.$

We now turn to the problem of relating the function $N(\sigma)$ to the
convex Szlenk index.

\begin{Lem}\label{Szest2} Suppose $0<\sigma<1.$  If $N= N(\sigma)$
 there exist $0<\epsilon_1,\ldots,\epsilon_{N+1}<1$ so that
$\sum_{k=1}^{N+1}\epsilon_k>\frac13\sigma^{-1}$ and
$\iota_{\epsilon_1}\iota_{\epsilon_2}\ldots\iota_{\epsilon_{N+1}}B_{X^*}\neq
\emptyset.$\end{Lem}

\begin{proof} Suppose $(x_a)_{a\in S}$ is a weakly null tree-map  of
height $N+1$ with
$\|x_a\|\le \sigma$ and so that $\|\sum_{a\in\beta}x_a\|>1$ for every
branch $\beta.$

Fix $\delta>0$ so that $(2N+3)\delta<\frac13.$ For any $a\in S$ with
$|a|=N+1$ we choose
$y_a^*$ with $\|y_a^*\|=1$ and $\langle \sum_{b\le a}x_b,y_a^*\rangle
=\|\sum_{b\le a}x_b\|.$  If $|a|>N$ set $y_a^*=y_b^*$ where $|b|=N$ and
$b\le a$. We now define
$y_a^*$ by backwards
induction so that for each $a$, $y_a^*$ is a weak$^*$-cluster point of
$\{y_b^*:b\in a+\}.$ It is then easy to apply induction to produce a full
tree $T\subset S$ so that $\lim_{b\in a+}y_b^*=y_a^*.$

Now let $x_a^*=y_a^*-y_{a-}^*$ when $|a|\ge 1$, and $x_{\emptyset}^*=0$ so
that
$(x_a^*)$ is a
weak$^*$-null tree-map in $X^*$ of height $N.$  Let
$y^*=y_{\emptyset}^*$
so that $\|y^*\|\le 1.$  We thus have on every branch $\beta$ of $T$
that
$$ \|y^*+\sum_{a\in\beta}x_a^*\|=1$$
and so
$$ \|\sum_{a\in\beta}x_a^*\|\le 2.$$

Now using Lemma \ref{reduce} we can pass to a further full tree
$T_1\subset T$ so that for suitable $\epsilon_1,\ldots,\epsilon_{N+1}$
we have
$2\epsilon_k\le \|x_a^*\| \le 2\epsilon_k+\delta$ if $|a|=k$ for $1\le
k\le N.$ Then by Proposition \ref{Szlenk} we have that
$\iota_{\epsilon_1}\ldots\iota_{\epsilon_{N+1}}B_{X^*}\neq\emptyset.$

Next it is clear that we can pass to a further full subtree $T_2$ so that
for every $a$ we have $|y^*(x_a)|<\delta$ (since $(x_a)$ is weakly null).
Finally we use Lemma \ref{dual} to produce a full tree $T_3\subset T_2$
so that for any branch $\beta$
$$ |\langle \sum_{a\in\beta}x_a,\sum_{a\in\beta}x^*_a\rangle
-\sum_{a\in\beta}\langle x_a,x_a^*\rangle|\le \delta.$$

Now for any branch $\beta$ in $T_3$ we have
$$ \|\sum_{a\in\beta}x_a\| =\langle\sum_{a\in\beta}x_a,
y^*+\sum_{a\in\beta}x_a^*\rangle.$$
Hence
\begin{equation}\begin{align*}1< \|\sum_{a\in\beta}x_a\| &\le
(N+2)\delta
+\sum_{a\in\beta}\langle
x_a,x_a^*\rangle\\  &\le (2N+3)\delta
+2\sigma\sum_{k=1}^{N+1}\epsilon_k.\end{align*}\end{equation}
Thus we have $\sum_{k=1}^{N+1}\epsilon_k\ge \frac13\sigma^{-1}.$

\end{proof}

\begin{Thm}\label{equiv}Suppose $X$ is a separable Banach space
containing no copy of
$\ell_1.$  Let $H(\tau)=(\text{Cz}(X,\tau)-1)^{-1}$ for $0\le \tau
<1.$  Then there is a universal constant $C\le 19200$ so that
$N(\sigma)^{-1}$ is $C$-equivalent to $\varphi(\sigma)$ and $H(\tau)$ is
$C$-equivalent to $\psi(\tau).$

\end{Thm}

\begin{proof} Suppose $0<\sigma<1$.  Then by Lemma \ref{summable} and
the preceding Lemma \ref{Szest2} we have that if
$\tau\le \frac1{12}\sigma^{-1}(N+1)^{-1}$ then
$$ \text{Cz}(X,\tau)-1 \ge \frac{1}{24\sigma\tau}.$$

Thus
$$H(\tau)\le 24\sigma\tau.$$  Then for any $\sigma$ if $0<
\tau \le \frac{25}{12}\sigma^{-1}(N(\sigma/25)+1)^{-1}$ we have
$$ H^*(\sigma) \ge \sigma\tau-
\frac{24}{25}\sigma\tau=\frac{1}{25}\sigma\tau.$$

Hence
$$ H^*(\sigma)\ge \frac1{12} (N(\sigma/25)+1)^{-1}\ge
\frac{1}{24}(N(\sigma/25))^{-1}.$$
Since $H^*$ is convex this implies that $H^*$ 600-dominates $N^{-1}.$

Now by Proposition \ref{orlicz} we have that $H$ 4-dominates $\psi$ and
hence $\psi^*$ 4-dominates $H^*$.  Thus $\varphi$ C-dominates
$N^{-1}$ with $C\leq 19200$.  By the remarks after Theorem \ref{UKK2} this
means that
$\varphi$ and $N^{-1}$ are C-equivalent.

Now recall that $H(n\tau)\ge nH(\tau)$ for $n\in\Bbb N$ by Lemma
\ref{convex}.

Thus for $s\ge 1$ we have $H(st)\ge \frac12sH(t).$ It follows easily that
$H$ is
2-equivalent to a function
$H_1$ with the
property that $H_1(t)/t$ is increasing, which is then 2-equivalent to a
convex function.  Hence $H$ is 4-equivalent to a convex function and
hence also to $H^{**}.$

Now $H^{*}$ 600-dominates $N^{-1}$ and hence 1200-dominates $\varphi.$
Thus $H^{**}$ is 1200-dominated by $\varphi^*$ and thus $4800-$dominated by
$\psi.$  Hence $H$ is 19200-dominated by $\psi$ and so $H$ is
19200-equivalent to $\psi.$
\end{proof}

\begin{Thm}\label{CzSz}Suppose $X$ is a separable Banach space not
containing $\ell_1.$  Then there is a universal constant $C<10^6$ so that if
$0<\tau\le 1,$
$$ \text{Cz}(X,\tau) \le \sum_{\substack{k\ge 0\\ 2^k\tau/C\le
1}}2^k\text{Sz}(X,2^k\tau/C).$$\end{Thm}

\begin{proof} Let $$K(\tau)=\left(\sum_{\substack{k\ge 0\\2^k\tau\le 1}}
2^k(\text{Sz}(X,2^k\tau)-1)\right)^{-1}.$$  Then arguing as with $H$ we
have that $K$ is 4-equivalent to $K^{**}.$

We next note that by Lemma \ref{Szest2} and Lemma \ref{summable} if
$0<\sigma<1$ we have that if $0<\tau\le
\frac16\sigma^{-1}(N(\sigma)+1)^{-1}$ then
$$ K(\tau)\le 6\sigma\tau.$$  Reasoning as above gives
$$  K^*(\sigma)\ge \frac16(N(\sigma/7)+1)^{-1}\ge
\frac1{12}N(\sigma/7)^{-1}.$$

Hence $K^*$ 100-dominates $N^{-1}$.  By Theorem \ref{equiv} this implies
that $K$ is $C$-dominated by $H$ for a suitable absolute constant $C\leq
600.100.4.4<10^6.$
This clearly yields the result.\end{proof}

Our next Corollary follows easily from Theorem \ref{CzSz}.

\begin{Cor}\label{power} Suppose $X$ is a separable Banach space with
$\text{Sz}(X)=\omega_0.$  Suppose  $p>1$ and
$$\sup_{0<\epsilon<1}\epsilon^p\text{Sz}(X,\epsilon)<\infty.$$ Then
$$ \sup_{0<\epsilon<1}\epsilon^p\text{Cz}(X,\epsilon)<\infty.$$\end{Cor}

We are now in position to state our main renorming theorems. The first one
simply
restates Theorem \ref{equiv}
:

\begin{Thm}\label{renorm} Suppose $X$ is a separable Banach space with
$\text{Sz}(X)=\omega_0.$   Then there exists an absolute constant $C\leq
19200$ such
that
 for any $0<\tau<1$ there is 2-equivalent norm $|\cdot|$ on $X$
so that if $x^*,x_n^*\in X^*$ satisfy $|x^*|=1$, $|x_n^*|=\tau$ and
$\lim_{n\to\infty}x_n^*=0$ weak$^*$ then
$$ \liminf_{n\to\infty}|x^*+x_n^*| \ge 1+
\frac{1}{\text{Cz}(X,\tau/C)}.$$

\end{Thm}

It is clear that Theorem \ref{renorm} is in a sense a best possible
result. Notice that we do not have a simultaneous renorming which works
for all $0<\tau<1.$   However we can combine these norms to give a single
norm which has the $UKK^*$ property with a nearly optimal modulus. Note
that it is clear that any space $X$
which has an equivalent $UKK^*$ norm satisfies $\text{Sz}(X)=\omega_0.$

\begin{Thm}\label{renorm2} Suppose $X$ is a separable Banach space with
$\text{Sz}(X)=\omega_0.$  Let $p_X$ be the Szlenk power type of $X.$  Then
for any $q>p_X$ there is an equivalent norm $|\cdot|$ on $X$ and a
constant $c>0$ so that if
$0<\tau\le 1$ and
$x^*,x_n^*\in X^*$ satisfy $|x^*|=1,$ $|x_n^*|=\tau$ and
$\lim_{n\to\infty}x_n^*=0$ weak$^*$ then
$$ \liminf_{n\to\infty}|x^*+x_n^*| \ge 1 +c\tau^q.$$\end{Thm}

\begin{proof} Fix $p_X<r<q.$ In this case we have that
$\text{Sz}(X,\epsilon)\le
C\epsilon^{-r}$ for some constant $C.$  By Corollary \ref{power} we have
a similar estimate $\text{Cz}(X,\epsilon)\le C\epsilon^{-r}.$  Hence, for
a suitable constant $c_1>0,$ for each $k\in\Bbb N$ there is a norm
$|\cdot|_k$ on
$X$ which is
2-equivalent to the original norm and such that if $|x^*|_k=1$,
$|x_n^*|_k=2^{-k}$ with $x_n^*$ weak$^*$-null then
$$ \liminf_{n\to\infty}|x^*+x_n^*|_k\ge 1+ c_12^{-rk}.$$

Now define the (dual) norm $|\cdot|$ on $X^*$ by
$$ |x^*|=\sum_{k=1}^{\infty}2^{(r-q)k}|x^*|_k.$$
This clearly defines an equivalent dual norm on $X^*$. Thus there is a
uniform constant $B$ so that for every $k$ we have $B^{-1}|x^*|\le
|x^*|_k\le B|x^*|.$ Suppose
$|x^*|=1$ and
$|x_n^*|=\tau$ with
$(x_n^*)$ weak$^*$-null.   Pick $k\in\Bbb N$ so that $2^{-k}\le
B^{-2}\tau\le 2^{1-k}.$  Then $|x^*|_k\le B$ and $|x_n^*|_k\ge
B^{-1}\tau\ge 2^{-k}|x^*|_k.$  Hence
$$ \liminf |x^*+x_n^*|_k\ge |x^*|_k(1 + c_12^{-rk}).
$$

This implies that
$$ \liminf |x^*+x_n^*| \ge 1+ c_2 2^{-qk} \ge 1+c_3\tau^q$$
for  suitable $c_2,c_3>0.$

\end{proof}

In general we do not know whether the functions $\text{Cz}(X,\tau)$
and $\text{Sz}(X,\tau)$ are equivalent.  However in certain cases, which
include for instance super-reflexive spaces but also James' quasi-reflexive
space $J$ (\cite{L3}),
they are equivalent:

\begin{Thm}\label{superr} Suppose $X$ is a separable Banach space with
$\text{Sz}(X)=\text{Sz}(X^*)=\omega_0.$  Then there
is a constant
$C$
(depending only on $X$) so that $\text{Cz}(X,\tau)\le
\text{Sz}(X,\tau/C).$\end{Thm}

\begin{proof} We begin by noting that $X^{**}$ is separable.  From
Theorem \ref{renorm}  and the fact that $Cz(X^*,\tau)<\infty$ for all
$\tau>0$ it follows that we can replace the original norm  with
an equivalent norm on
$X\subset
X^{**}$ so that there exists $\delta>0$ with the property that if
$\|x\|=\|x_n\|=1$ and
$(x_n)$ is weakly null then
\begin{equation}\label{prop} \liminf \|x+x_n\|>
1+\delta.\end{equation} Using this property we make an estimate of
$N(\sigma).$ If
$N=N(\sigma)<\infty$ there is weakly null tree-map $(x_a)_{a\in\cal
F\mathbb N}$ of height
$N+1$ so that $\|x_a\|\le \sigma$ for every $a\in S$ but
$\|\sum_{a\in\beta}x_a\|>1$ on every branch.  We define a second tree-map
$(y_a)_{a\in\cal F\mathbb N}$ by
$y_a=x_a$ if $|a|\le N+1$ and $y_{a}=x_{b}$ if $N+1<|a|=k\le 2(N+1)$ and
$b=\{m_{N+2},\ldots,m_k\}$ if $ a=\{m_1,\ldots,m_k\}.$  If $|a|>2(N+1)$
then $y_{a}=0.$  It follows easily from (\ref{prop}) that there is a full
subtree $S$ so that on every branch we have $\sum_{a\in\beta}>1+\delta.$
In particular $N((1+\delta)^{-1}\sigma)\le 2N(\sigma)+1\le 3N(\sigma).$

Iterating
we have $N(\sigma)<\infty$ for all $\sigma$.
The estimate on $N$ clearly implies an estimate of the type
$$ N(\lambda\sigma) \le C\lambda^{-p}N(\sigma)\qquad
0<\lambda,\sigma<1.$$ where $1< p<\infty.$

By Theorem \ref{equiv} and Proposition \ref{orlicz} we obtain a dual
estimate $\text{Cz}(X,\lambda\tau)\ge C\lambda^{-q}\text{Cz}(X,\tau)$
where $\frac1p+\frac1q=1.$ Note that $q>1.$

We now use Theorem \ref{CzSz}.  First note that for a suitable constant
$C_1$ we have
\begin{equation}\begin{align*}
\sum_{\substack{k\ge 0\\ 2^k\tau\le 1}}2^k\text{Sz}(X,2^k\tau) &\le
 \sum_{\substack{k\ge 0\\ 2^k\tau\le 1}}2^k\text{Cz}(X,2^k\tau) \\
&\le
C \sum_{\substack{k\ge 0\\ 2^k\tau\le 1}}2^{k(1-q)}\text{Cz}(X,\tau) \\
&\le  C_1\text{Cz}(X,\tau).\end{align*}\end{equation}

Recall we also have an estimate from Theorem \ref{CzSz}
$$ \text{Cz}(X,\tau)\le \sum_{\substack{k\ge0\\ 2^k\tau/C\le 1}} 2^k\text
{Sz}(X, 2^k\tau/C_2).$$ Now pick $k_0$ so that for every $\tau$ we have
$$  \text{Cz}(X,2^{-k_0}C_2\tau) \ge (1+2^{k_0}C_1)\text{Cz}(X,\tau).$$

This is possible by the growth condition on $\text{Cz}(X,\tau)$ as
$\tau\to 0.$ Then
\begin{equation}\begin{align*}
\text{Cz}(X,\tau) &\le
\text{Cz}(X,2^{-k_0}C_2\tau)-2^{k_0}C_1\text{Cz}(X,\tau)\\
&\le  \sum_{\substack{k\ge 0\\2^{k-k_0}\tau\le 1}}
2^k\text{Sz}(X,2^{k-k_0}\tau)  -\sum_{\substack{k\ge 0\\ 2^k\tau\le 1}}
2^{k+k_0}\text{Sz}(X,2^k\tau) \\
&\le \sum_{k=0}^{k_0-1}2^k\text{Sz}(X,2^{k-k_0}\tau)\\
&\le 2^{k_0}\text{Sz}(X,2^{-k_0}\tau).\end{align*}\end{equation}

This implies the Theorem.\end{proof}

Note that the above proof shows that if $X$ is separable and
$\text{Sz}(X)=\text{Sz}(X^*)=\omega_0$ then $N(\sigma)<\infty$ for any
$\sigma>0$. Hence by Theorem \ref{sumszlenk} below, the Szlenk index of
such spaces is not summable.

We conclude this section with a characterization of spaces with summable
Szlenk index (see also Theorem \ref{UKK2}). Recall that by Lemma
\ref{summable}, condition (iii) below is equivalent to the
summability of the convex Szlenk index Cz. Note that Sz and Cz are also
equivalent for the spaces which satisfies the conditions of
the next theorem.

\begin{Thm}\label{sumszlenk}  Let $X$ be a separable Banach space.  The
following assertions are equivalent:
\newline
(i) $X$ has summable Szlenk index.\newline
(ii) There exists $\sigma>0$ so that $N(\sigma)=\infty.$\newline
(iii) There exists $K>0$ so that for any $0<\tau<1$,
$\text{Cz}(X,\tau)\le K\tau^{-1}.$\newline
(iv) For any function $f:(0,1)\to (0,1)$ which satisfies $\lim_{\tau\to
0}\tau^{-1}f(\tau)=0,$ there is an equivalent norm
$|\cdot|$ on
$X$ and a constant $c>0$ so that if $0<\tau<1$ and $x^*,x_n^*\in X$
satisfy $|x^*|=1$, $|x_n^*|=\tau$ and $\lim_{n\to\infty}x_n^*=0$ weak$^*$
then $$\liminf_{n\to\infty}|x^*+x_n^*|\ge 1+cf(\tau).$$

\end{Thm}

\begin{proof} We first note that (iii) $\Rightarrow$ (i) follows directly
from the second part of Lemma \ref{summable}.

Next we prove that (iv) $\Rightarrow$ (iii).  Assume that (iii) does not
hold.  Let $$ f(\tau)=\sqrt{\frac{\tau}{\text{Cz}(X,\tau)-1}}.$$  Then we
have $\lim_{\tau\to 0}\tau^{-1}f(\tau)=0.$  Thus assuming (iv) we
conclude an estimate $\text{Cz}(X,\tau)\le C(f(\tau))^{-1}$ which gives a
contradiction.

Next we note that (i) $\Rightarrow$ (ii) is immediate from Lemma
\ref{Szest2}.

It remains therefore to show that (ii) $\Rightarrow$ (iv).  Note first
that we can assume that $\tau^{-1}f(\tau)$ is a monotone increasing
function. It will be sufficient to prove the result with $f(\tau)$
replaced by $f(\tau/4).$

We note that by Theorem \ref{equiv} and Proposition \ref{orlicz} we have
that $\psi(\tau)\ge c\tau$ for some $c>0.$  Hence if $0<c_1<c$ there is,
for each $k\in\mathbb N,$ a 2-renorming $|\cdot|_k$ of $X$ so that
$\frac12\|x\|\le |x|_k\le \|x\|$ for $x\in X$ and
$$ \liminf|x^*+x_n^*|_k \ge 1+c_12^{-k}$$
whenever $|x^*|_k=1$ and $(x_n^*)$ is a weak$^*$-null sequence in $X^*$
with
$|x^*_n|_k=2^{-k}.$  Now using the same idea as in the preceding proof we
define a dual norm on $X^*$ by
$$
|x^*|=\frac{1}{2f(\frac12)}\sum_{k=1}^{\infty}\left(2^kf(2^{-k})-
2^{k+1}f(2^{-k-1})\right)
|x^*|_k.$$

It follows from our definition that $|\cdot|$ is well-defined and
2-equivalent to the original norm. Now suppose $x^*,x_n^*\in X^*$
are such that
$|x^*|=1$, $|x_n^*|=\tau$ and $(x_n^*)$ is weak$^*$-null.  Pick $j\in
\mathbb N$ so that $2^{-j}\le \frac14\tau\le 2^{1-j}.$  Then if $k\ge j$
we have $|x^*|_k \le 2$ and $|x_n^*|_k\ge \frac12\tau\ge 2^{-j}|x^*|_k.$

Hence we obtain by convexity
$$ \liminf |x^*+x_n^*|_k \ge |x^*|_k(1+c_12^{-j}).$$
Summing over $k$ we have
$$ \liminf |x^*+x_n^*| \ge 1 + c_12^{-j}\sum_{k=j}^{\infty}
(2^kf(2^{-k}) - 2^{k+1}f(2^{-k-1}))|x^*|_k.$$

Hence for a suitable constant $c_2>0$ we have
$$ \liminf |x^*+x_n^*|\ge 1 + c_2f(\tau/4).$$
This completes the proof.\end{proof}

{\it Remarks.} Spaces with summable Szlenk index are considered in
\cite{KOS} where it is shown that the original Tsirelson space is a
reflexive example.  These spaces are very close to being subspaces of
$c_0$ which are characterized by the existence of an equivalent
Lipschitz-UKK$^*$-norm (cf. \cite {KW}, \cite{GKL}) i.e. a norm such that
$$ \liminf \|x^*+x_n^*\| \ge 1+ c\tau$$
whenever $\|x^*\|=1$ and $(x_n^*)$ is a weak$^*$-null sequence with
$\|x_n^*\|=\tau.$ The ``lack of isotropy" of the Szlenk derivation
when applied to the Tsirelson space seems to be responsible for the difference
between summability of the Szlenk index and embeddability into $c_0$.

\section{Applications to uniform and Lipschitz homeomorphisms}

We first recall the Gorelik principle.  We use the version proved in
\cite{GKL}; see \cite{JLS} for an earlier version. Both are based on the
original idea of Gorelik \cite{G}. If $V:X\to Y$ is uniformly continuous
we denote the modulus of continuity by $\omega(V,t)=\sup \{\|Vx-Vy\|;\
\|x-y\|\le t\}.$

\begin{Prop}\label{Gorelik} (The Gorelik Principle) Let $X$ and $Y$ be
two Banach spaces and let $U$ be a homeomorphism of $X$ onto $Y$ whose
inverse is uniformly continuous.  Suppose $b>0$ and $d>\omega(U^{-1},b)$
and let $X_0$ be a closed subspace of $X$ of finite codimension.  Then
there is a compact subset $K$ of $Y$ so that
$$ bB_Y \subset K+ U(2d B_{X_0}).$$
\end{Prop}

Let us recall at this point that if $V:X\to Y$ is uniformly continuous it
is Lipschitz at large distances.  More precisely if $0<t<\infty$ then
we have
$$ \|Vx-Vy\|\le t^{-1} \omega(V,t) (t+\|x-y\|).$$
Let us define the asymptotic Lipschitz constant of $V$ by
$$ \|V\|_u =\lim_{t\to\infty}t^{-1}\omega(V,t).$$
If $X$ and $Y$ are uniformly homeomorphic Banach spaces we define
$$d_u(X,Y)=\inf\{\|U\|_u\|U^{-1}\|_u:\ U:X\to Y \text{ is a uniform
homeomorhism}\}.$$

Its clear that $d_u(X,Y)\ge 1$ and that $d_u(X,Z)\le d_u(X,Y)d_u(Y,Z)$ if
$X$, $Y$ and $Z$ are all uniformly homeomorphic.

\begin{Lem}\label{scaling} Suppose $X,Y$ are uniformly homeomorphic.
Then if $M^2>d_u(X,Y)$  and $\eta>0$ there is a uniform
homeomorphism
$V:X\to Y$ so that we have
\begin{equation}\begin{align*}
 \|Vx_1-Vx_2\| &\le M\max(\|x_1-x_2\|,\eta) \qquad & x_1,x_2\in X\\
\|V^{-1}y_1-V^{-1}y_2\| &\le M\max(\|y_1-y_2\|,\eta) \qquad &y_1,y_2\in
Y.
\end{align*}\end{equation}

\end{Lem}

\begin{proof} Let $U:X\to Y$ be a uniform homeomorphism with
$\|U\|_u\|U^{-1}\|_u=L^2<M^2.$  For any $s>0$ define
$$ V_s(x)= L\|U\|_u^{-1}s^{-1}U(sx)$$
so that
$$ V_s^{-1}(y)= s^{-1} U^{-1}(L^{-1}\|U\|_u sy).$$
Thus $\omega(V_s,t) = L\|U\|_u^{-1}s^{-1}\omega(U,st)$ and
$\omega(V_s^{-1},t)=s^{-1}\omega(U^{-1},L^{-1}\|U\|_ust).$   Hence for any
$\delta>0$ we have
$$ \|V_sx_1-V_sx_2\| \le L\|U\|_u^{-1}\delta^{-1}s^{-1}\omega(U,s\delta)
(\delta+\|x_1-x_2\|) $$
and
$$ \|V_s^{-1}y_1-V_s^{-1}y_2\| \le
\delta^{-1}s^{-1}\omega(U^{-1},L^{-1}\|U\|_us\delta)(\|y_1-y_2\|+\delta).$$
It is clear that if we take $s$ large enough and $\delta$ small enough,
we obtain the lemma.\end{proof}

We now state our main result of this section:

\begin{Thm}\label{uniform} Suppose $X$ and $Y$  are separable Banach
spaces which are
uniformly  homeomorphic.    Then for any $\lambda>d_u(X,Y)$
and any $\epsilon>0$ there is a
$\lambda$-equivalent norm
$|\cdot|$ on $Y$ with $$\theta_{Y,|\cdot|}(\tau)\ge
\theta_X(\tau/C)-\epsilon,\qquad \forall 0\le \tau\le 1,$$ where
$C=32\lambda^2.$

\end{Thm}

 \begin{proof} Let $M=\lambda^{1/2}.$ We assume that $U:X\to Y$ is a
uniform homeomorphism of
$X$ onto $Y$ which satisfies:
\begin{equation}\begin{align*}
 \|Ux_1-Ux_2\| &\le M\max(\|x_1-x_2\|,1) \qquad & x_1,x_2\in X\\
\|U^{-1}y_1-U^{-1}y_2\| &\le M\max(\|y_1-y_2\|,1) \qquad &y_1,y_2\in Y.
\end{align*}\end{equation}

This is possible by Lemma \ref{scaling}.  Recall that
$C=32\lambda^2=32M^4.$

We first define a decreasing sequence of dual norms $\{|\cdot|_k
\}_{k=1}^{\infty}$ on $Y^*$ by
$$ |y^*|_k =\sup \left\{ \frac{|y^*(Ux_1-Ux_2)|}{\|x_1-x_2\|};\quad
x_1,x_2\in X,\ \|x_1-x_2\|\ge 2^k\right\}.$$
It is clear that we have $M^{-1}\|y^*\|\le |y^*|_k \le M\|y^*\|$
for
every $k\ge 1.$

We will prove the following Claim concerning the norms $|\cdot|_k$:

{\it Claim:  suppose $0<\tau<1$ and  $y^*,y_n^*\in Y^*$ are such that
$\|y^*\|\le M$, $(y_n^*)$ is weak$^*$-null and $\|y_n^*\|\ge M^{-1}\tau$
for all $n.$  Then  if $2^kC\tau^{-1}\theta_X(C^{-1}\tau)> 4M,$
\begin{equation}\label{claim}
\liminf_{n\to\infty}|y^*+y_n^*|_k \ge
2|y^*|_{k+1}-|y^*|_k+\theta_X(\frac{
\tau}{C}).\end{equation} }

To prove the claim we first define $\sigma=C\tau^{-1}\theta_X(\tau/C)$ so
that $0<\sigma\le 1.$ Let us also for convenience of exposition write
$\beta=\theta_X(\tau/C).$ Now for any
$\epsilon>0$ we can choose
$x,x'\in X$ so that $\|x-x'\|\ge 2^{k+1}$ and $y^*(Ux-Ux')\ge
(1-\epsilon)\|x-x'\||y^*|_{k+1}.$  As usual in such arguments we can
suppose by using translations that $x'=-x$ and $Ux'=-Ux.$ Next we apply
Proposition \ref{renorming} and the separability of $X$.  We deduce
the existence of a finite codimensional subspace $X_0$ of $X$ so that
\begin{equation}\label{1}
\|x+z\|\ge \|x\|\ge 2^k \qquad \forall z\in X_0\end{equation} and
\begin{equation}\label{2} \|x+z\|\le (1+2\beta)\|x\|\qquad
\forall z\in\sigma\|x\|B_{X_0}. \end{equation}

Now, by assumption we have $\omega(U^{-1},b)< 2Mb$ if $b\ge 1.$  Hence we
apply the Gorelik principle, Proposition \ref{Gorelik} for
$b=\frac{\sigma\|x\|}{4M}>1$ by choice of $k,$ and
$d=\frac{\sigma\|x\|}{2}.$
We deduce that there is a compact subset $K$ of $Y$ so that:
\begin{equation}\label{3}
\frac{\sigma \|x\|}{4M}B_Y \subset K+
U(\sigma\|x\|B_{X_0}).\end{equation}

It follows from (\ref{3}) that there is a sequence $z_n\in
\sigma\|x\|B_{X_0}$ so that
$$ \liminf_{n\to\infty}y_n^*(-Uz_n) \ge \frac{\sigma\tau}{4M^2}\|x\|=
\frac{C\beta}{4M^2}\|x\|.$$

Now $\|x+z_n\|=\|z_n-x'\|$ so that by (\ref{1}) and (\ref{2}) we have
$$ y^*(Ux+Uz_n)=y^*(Uz_n-Ux')\le
(1+2\beta)|y^*|_k\|x\|$$

But

$$ y^*(Ux)=\frac12y^*(Ux-Ux') \ge (1-\epsilon)|y^*|_{k+1}\|x\|.$$

Hence we deduce
$$ y^*(Uz_n) \le
\left[(1+2\beta)|y^*|_k-(1-\epsilon)
|y^*|_{k+1}\right]\|x\|.$$

Combining these estimates gives:
$$ \liminf_{n\to\infty}(y^*+y_n^*)(Ux-Uz_n)\ge
(2-2\epsilon)|y^*|_{k+1}\|x\|-|y^*|_{k}\|x\| +
\frac{C\beta}{4M^2}\|x\|-2\beta|y^*|_k\|x\|.$$

The left-hand side is estimated by
$(1+2\beta)\liminf|y^*+y_n^*|_{k}\|x^*\|.$  Hence since
$\epsilon>0$ is arbitrary and $\|y^*\|\le M,$
\begin{equation}\begin{align*}
 \liminf_{n\to\infty}|y^*+y_n^*|_k &\ge
(1-2\beta)(2|y^*|_{k+1}-|y^*|_k+\frac{C-8M^3}{4M^2}\beta) \\
&\ge 2|y^*|_{k+1}-|y^*|_k +
\frac{C-8M^3-16M^4}{4M^2}\beta.\end{align*}\end{equation}
Our choice of $C$ gives
$$ \frac{C-8M^3-16M^4}{4M^2}\ge 1$$
and so the Claim (\ref{claim}) is proved.

The completion of the argument from the claim is easy.  If
$\theta_X(1/C)\le \epsilon$ the original norm will suffice.  Otherwise
choose $\tau_0$ so that $\theta_X(\tau_0/C)=\epsilon.$

 Pick $k_0\in\mathbb N$ so that
$\frac{C}{\tau_0}\theta_X(\frac{\tau_0}{C})2^{k_0}>4M.$  Pick an integer
$N$ so that
$$ \frac{2M^2}{N}<\epsilon.$$

Now let $$|y^*|=\frac1N \sum_{k=k_0+1}^{k=k_0+N}|y^*|_k$$ which clearly
defines a dual norm on $Y^*$ with $M^{-1}\|y^*\|\le |y^*|\le M\|y^*\|.$

Assume $|y^*|=1$ and that $(y_n^*)$ is weak$^*$ null with $|y_n^*|\ge
\tau.$  Then $\|y^*\|\le M$ and $\|y_n^*\|\ge \tau/M.$  Hence for each
$k_0+1\le k\le k_0+N$ we have (\ref{claim}).  Summing gives
$$ \liminf |y^*+y_n^*|\ge |y^*|
-\frac2N|y^*|_{k_0+1}+\theta_X(\tau/C).$$

Now
$$ \frac2N|y^*|_{k_0+1}\le \frac{2M}{N}\|y^*\|\le
\frac{2M^2}{N}<\epsilon$$

so that the theorem is proved.\end{proof}

Before turning to applications for this result, let us notice that in
the
case of Lipschitz equivalence the same techniques can be applied to give
a rather stronger result. A special case of the theorem below is shown in
\cite{GKL}.

\begin{Thm}\label{lipschitz}  Suppose $X$ and $Y$ are separable
Banach
spaces which are Lipschitz isomorphic.  Then there is an equivalent norm
on $Y$ so that $\theta_Y$ $C$-dominates $\theta_X$ for some
$C>0.$\end{Thm}

\begin{proof}  We only provide a sketch.  In this case one can use the
norm on $Y^*$ defined by
$$
|y^*|=\sup\left\{\frac{|y^*(Ux_1)-y^*(Ux_2)|}{\|x_1-x_2\|}\qquad
x_1\neq x_2\right\}.$$

The calculations are similar but rather simpler.  We leave the details to
the reader.\end{proof}

Our main application of Theorem \ref{uniform} is the fact that the convex
Szlenk index is (up to equivalence) invariant under uniform
homeomorphisms.

\begin{Thm}\label{czinv}  Suppose $X$ and $Y$ are uniformly homeomorphic.
Then:\newline
(i) $\text{Sz}(X)\le \omega_0$ if and only if $\text{Sz}(Y)\le
\omega_0$.\newline
(ii)  There exists a constant $C$ so that if $0\le \tau\le 1$
$$ \text{Cz}(X,C\tau)\le \text{Cz}(Y,\tau)\le \text{Cz}(X,\tau/C).$$
\end{Thm}

\begin{proof}  We note first that the relevant functionals are separable
determined (\cite{L2}). Hence we may assume that $X$
and $Y$ are separable. The result is then an almost immediate deduction
from Theorem
  \ref{uniform}.    In fact it is immediately clear that the statement of
the Theorem yields that for the original norm on $Y,$
$$ \text{Cz}(Y,\tau) \le   (\theta_X(\frac{\tau}{C\lambda}))^{-1}+1.$$
 If we then consider all $2$-equivalent norms on $X$ we deduce an
estimate of the form
$$ \text{Cz}(Y,\tau) \le \psi_X(\tau/C)^{-1}+1$$
for a different constant $C.$  However by Theorem \ref{renorm} this
implies an estimate
$$ \text{Cz}(Y,\tau)\le \text{Cz}(X,\tau/C).$$   This estimate and its
converse establish both (i) and (ii).\end{proof}

\begin{Thm}\label{c0}(i) Let $X$ be a Banach space which is uniformly
homeomorphic
to a subspace of $c_0$.  Then $X$ has summable Szlenk index (and in
particular, $X^*$ is separable). \newline
(ii) If $X$ is uniformly isomorphic to $c_0$ then $X^*$ is linearly
isomorphic to $\ell_1.$\end{Thm}

\begin{proof}  Since any subspace of $c_0$ has summable Szlenk index then
$X$ must also have summable Szlenk index by Theorem \ref{czinv} and
Theorem \ref{sumszlenk}.  For the second part note that $X$ must be a
$\cal L_{\infty}-$space \cite{HM} and so (i) yields that $X^*$ is
isomorphic to $\ell_1$ \cite{LS}. \end{proof}

{\it Remarks.}1) It is shown in \cite{A} that the Bourgain-Delbaen preduals
of $l^1$ have a
finite Szlenk index, although they do not contain $c_0$. On the other hand,
it is shown in
\cite{H} that these spaces contain hereditarily $l^p$ for some $p\in
(1,\infty)$. It follows now
from Theorem \ref{c0} that no infinite-dimensional subspace of a
Bourgain-Delbaen space is
uniformly homeomorphic to a subspace of $c_0$.

2)As remarked above, there are unfortunately examples of
Banach spaces
with summable Szlenk index which do not embed into $c_0$ as shown in
\cite{KOS}.  We conjecture however that any predual of $\ell_1$ with
summable Szlenk index is isomorphic to $c_0.$

A Banach space which is Lipschitz isomorphic to $c_0$ is (in the separable
case) linearly isomorphic to $c_0$ (\cite {GKL}). We do
not know whether a space which is uniformly homeomorphic to
$c_0$ is isomorphic to
$c_0$. Curiously however, we can show that if
$d_u(X,c_0)$ is small enough then
$X$ is isomorphic to $c_0.$ This result is very similar to the
quantitative Lipschitz result in \cite {GKL}. Note that the following
statement is true for any equivalent renorming of
$c_0$, although the function $f$ and $\epsilon_0$ depend upon the given norm.

\begin{Thm}\label{close} There is $\epsilon_0>0$ and a function
$f:(0,\epsilon_0]\to (0,\infty)$ with $\lim_{\epsilon\to
0}f(\epsilon)=0$ so that if $d_u(X,c_0)<1+\epsilon$ then
$d(X,c_0)<1+f(\epsilon).$ \end{Thm}

\begin{proof} We use the ideas of \cite{KO}: in particular we use the
 the Gromov-Hausdorff distance for two Banach spaces $d_{GH}(X,Y),$ for
which we refer to \cite{KO}.  From Theorem 5.9 of \cite{KO} we have the
fact that if $X_n$ is a sequence of Banach spaces converging to $c_0$ in
$d_{GH}$ then the Banach-Mazur distance $d(X_n,c_0)\to 1.$  It therefore
suffices to show that $d_u(X_n,c_0)\to 1$ implies $d_{GH}(X,c_0)\to 0.$

To achieve this suppose that $d_u(X,Y)<1+\epsilon.$  Then for any
$\eta>0$ we use Lemma \ref{scaling} to produce a uniform homeomorphism
$V:X\to Y$ such that $V(0)=0$ and
\begin{equation}\begin{align*}
 \|Vx_1-Vx_2\| &\le (1+\epsilon)^{1/2})\max(\|x_1-x_2\|,\eta) \qquad
& x_1,x_2\in X\\
\|V^{-1}y_1-V^{-1}y_2\| &\le
(1+\epsilon)^{1/2})
\max(\|y_1-y_2\|,\eta)
\qquad
&y_1,y_2\in Y.
\end{align*}\end{equation}

 Now define $\Phi:B_X\to B_Y$ by $\Phi(x)=(1+\epsilon)^{-1/2}Vx$ and
$\Psi:B_Y\to B_X$ by $\Psi(y)=(1+\epsilon)^{-1/2}V^{-1}y.$

Then if $x\in B_X,\ y\in B_Y,$
$$
\left|\|x-\Psi(y)\|-\|x-V^{-1}(y)\|\right| \le \frac12\epsilon$$
and
$$\left |\|\Phi(x)-y\|-\|Vx-y\|\right|\le \frac12\epsilon.$$
Now
$$ \|Vx-y\| \le (1+\epsilon)^{1/2}(\|x-V^{-1}y\| +\eta)$$
and a similar reverse inequality gives
$$\left |\|Vx-y\|-\|x-V^{-1}y\|\right| \le
\frac12\epsilon(\max\|x-V^{-1}y\|,\|Vx-y\|) + 2\eta.$$

Since $\|x-V^{-1}y\|,\|Vx-y\|\le (1+\epsilon)^{1/2}+1<3$ we obtain
$$\left |\|x-\Psi(y)\|-\|\Phi(x)-y\|\right|<\frac52\epsilon+2\eta.$$
This implies (\cite{KO}) that $d_{GH}(X,Y)\le\frac54\epsilon+\eta.$
The result then follows.

\end{proof}

We now give one more application of our main result on the invariance of
the convex Szlenk index:

\begin{Thm}\label{ellp} Suppose $2<p<\infty$ and $X$ is a quotient (resp.
subspace) of $\ell_p.$  If
$Y$ is uniformly homeomorphic to $X$  then
$Y$ is linearly isomorphic to a quotient
(resp. subspace) of $\ell_p.$\end{Thm}

 {\it Remark.} In fact the methods of \cite{JLS} give the subspace
version rather easily so the quotient case is the more interesting.

\begin{proof} In both cases we have $\text{Cz}(X,\epsilon)\le
C\epsilon^{-q}$ where $\frac1p+\frac1q=1.$  We therefore deduce
by Theorem \ref{czinv} that $\text{Cz}(Y,\epsilon)\le C\epsilon^{-q}.$

Now we use the ``standard ultraproduct technique'' (cf.\cite {JLS} p. 438
or
\cite{B}).  The spaces $X$ and $Y$ are super-reflexive and hence we can
find Lipschitz-isomorphic separable spaces $X_1\supset X,Y_1\supset Y$ so
that
$X_1$
 is one-complemented in an ultraproduct $X_{\cal U}$ and $Y_1$ is
one-complemented in $Y_{\cal U}.$  Then $Y_1$ embeds complementably
into $X_1$ and $X_1$ embeds complementably in $Y_1.$ Since $X$ and $Y$
are complemented in their ultraproducts this leads to the fact that
 $Y$ is a quotient
(respectively a subspace) of $L_p$.

In the case when $X$ is a subspace of $\ell_p$ we can complete the
argument
very simply.  Since $p>2,$  if $Y$ does not embed into $\ell_p$ then
$\ell_2$ embeds complementably in $Y$ \cite{JO} and so
$\text{Cz}(X,\epsilon)\ge
c\epsilon^{-2}$ for some $c>0$ which yields a contradiction.

Now consider the case when $X$ is a quotient of $\ell_p$ so that $Y$
is a quotient of $L_p[0,1].$  We have
$\text{Cz}(Y,\epsilon)\le
C\epsilon^{-q}.$  It follows from Theorem \ref{equiv} and Proposition
\ref{orlicz} that $\psi_Y(\tau)\le C_1\epsilon^p$  for a suitable
constant $C_1.$  We will argue that this implies $Y$ is of type
$p$-Banach-Saks in the sense of \cite{J}, i.e. there is a constant
$\lambda>0$ so that every normalized weakly
null basic sequence $(w_n)_{n=1}^{\infty}$ has a subsequence $(v_n)$
satisfying $\|\sum_{i=1}^nv_i\| \le \lambda n^{1/p}.$

To do this we note that by the definition of $\psi_Y$ there is, for each
$k\in\mathbb N$ an equivalent norm $|\cdot|_k$ on $Y$ so that $\|y\|\le
|y|_k\le 2\|y\|$ for $y\in Y$ and $$\limsup |y+y_n|_k^p \le
|y|_k^p(1+2C_1k^{-1})$$ if $(y_n)$ is weakly null and $\sup|y_n|_k\le
k^{-1/p}|y|_k.$

Now if $(w_n)$ is a normalized weakly null basic sequence in
$(Y,\|\cdot\|)$,
we pass to a subsequence $(v_n)$ such that if $1\le k\le 2^n$ is such
that $|v_n|_k \le k^{-1/p}|\sum_{i=1}^{n-1}v_i|_k$ then
$$ |\sum_{i=1}^nv_i|_k \le |\sum_{i=1}^{n-1}v_i|_k(1+3C_1k^{-1})^{1/p}.$$

For fixed $n$, let $m=[\log_2n].$  Let $r$ be the greatest integer with
$1\le r\le n$ so
that $|v_{r}|_n> n^{-1/p}|\sum_{i=1}^{r-1}v_i|_n.$
 Provided
$1\ge
2n^{-1/p}(\log_2n+1),$ we have $r>m.$  Hence for $r\le j\le n-1$ we have
$$ |\sum_{i=1}^{j+1}v_i|_n \le (1+3C_1n^{-1})^{1/p}|\sum_{i=1}^jv_i|_n$$
and so \begin{equation}\begin{align*} \|\sum_{i=1}^nv_i\| &\le
 |\sum_{i=1}^n v_i|_n\\ &\le (1+3C_1n^{-1})^{n/p}|\sum_{i=1}^rv_i|_n \\
 &\le e^{3C_1/p}(|\sum_{i=1}^{r-1}v_{i}|_n+|v_r|_n)\\ &\le
e^{3C_1/p}(n^{1/p}+1)|v_r|_n \\
 &\le 4e^{3C_1/p}n^{1/p}.\end{align*}\end{equation}

We can now combine Theorems III.1 and III.2 (or Remark III.3) of \cite{J}
to deduce that
$Y$ is isomorphic to a quotient of $\ell_p.$\end{proof}

\end{document}